\documentclass[12pt]{amsart}

\newtheorem{thm}{Theorem}[section]
\newtheorem{prop}[thm]{Proposition}
\newtheorem{lem}[thm]{Lemma}
\newtheorem{cor}[thm]{Corollary}

\newtheorem{fact}[thm]{Fact}
\date{}

\def\figbox#1{\medskip
\begin{center} \fbox{#1} \end{center}
\medskip}

\newcommand{\fr}[1]{\mbox{${\mathfrak a}_{#1}$}}

\def\figbox#1{\medskip
\begin{center} \fbox{#1} \end{center}
\medskip}

\newcommand{\spf}{\noindent {\it Proof. }}

\newcommand{\epf}[1]{
\hfill \fboxsep = 1.5pt \fbox{{\footnotesize #1}}
\fboxsep = 3pt
\vspace{3mm}}
\newcommand{\fpf}{
\hfill \fboxsep = 1.5pt 
\fbox{{\footnotesize \thethm}}
\fboxsep = 3pt
\vspace{3mm}}

\begin{document}

\title{ALMOST ALTERNATING DIAGRAMS AND FIBERED 
            LINKS IN $S^{3}$}

\def\kyaku{The first author was supported in part by Hyogo 
Science and Technology Association. 
The second author is
supported by Fellowships
of the Japan Society for the Promotion of 
Science for Japanese Junior Scientists.}

\author{Hiroshi GODA, Mikami Hirasawa and Ryosuke YAMAMOTO}
\thanks{\kyaku}

\maketitle

\section{Introduction and statements of results}

The concept of Murasugi sum (for the definition, see Section 2)
of Seifert surfaces in the $3$-sphere $S^{3}$ was introduced by
K. Murasugi, and it has been playing important roles in
the studies of Seifert surfaces and links.
The Murasugi sum is known to be natural in many senses,
and in particular the following is known.
(We say that a Seifert surface $R$ is a {\it fiber surface} if
$\partial R$ is a fibered link and $R$ realizes the fiber.)

\begin{thm}[\mbox{\cite[Theorem 3.1]{G}}] \label{thm:fiber-sum}
Let $R$ be a
Murasugi sum of $R_{1}$ and $R_{2}$.
Then $R$ is a fiber surface if and only if
both $R_{1}$ and $R_{2}$ are fiber surfaces.
\end{thm}

On the other hand, the concept of alternating link has also
been important in knot theory.
It has been known that there are some relationships between
alternating diagrams and the Seifert surfaces obtained by
applying Seifert's algorithm to them.                    
For example, if a link diagram
$D$ is alternating, then the Seifert surface
obtained from $D$ by the algorithm is of minimal genus, 
\cite{C, Mu1958}.

In \cite{G}, D. Gabai gave a geometric proof to the following
theorem, which also follows
from \cite{Mu} and \cite{St}.
Note that if $L$ is fibered, then minimal genus Seifert surfaces
for $L$ are unique up to isotopy and the fiber is realized 
by the minimal genus surface.

\begin{thm}[\mbox{\cite[Theorem 5.1]{G}}] \label{thm:fiber-hopf}
Let $L$ be an oriented link with an alternating diagram $D$.
$L$ is a fibered link if and only if the surface $R$ obtained
by applying Seifert's algorithm to $D$ is connected and
(obviously) desums into a union of Hopf bands.
\end{thm}

We say that a Seifert surface $R (\subset S^3)$
{\it desums} into $R_{1},\ldots , R_{n}$ if
$R$ is a Murasugi sum of them. 
Especially, if $R$ is obtained by successively plumbing 
(i.e., 4-Murasugi summing)
finite number of Hopf bands to a disk, 
we call $R$ a {\it Hopf plumbing}.

Actually, the \lq only if' part of Theorem \ref{thm:fiber-hopf} 
can be strengthened as in the following theorem,
which follows from Propositions \ref{prop:para-deplumb}
and \ref{prop:alt-para}.

\begin{thm} \label{thm:alt-fiber}
Let $L$ be an oriented link with an alternating diagram $D$.
$L$ is a fibered link if and only if the surface $R$ obtained
by applying Seifert's algorithm to $D$ is a Hopf plumbing.
Moreover, $R$ is a fiber surface if and only if
$R$ is deformed into a disk by 
successively
cutting one of a pair of \lq parallel bands' 
(defined in Section 5).
\end{thm}

In \cite{A}, C. Adams et al. generalized the concept of alternating
links and introduced the concept of almost alternating links.
A diagram $D$ in $S^2$ is called
{\it almost alternating}
(resp. {\it $2$-almost alternating})
if $D$ becomes an alternating diagram after one crossing change
(resp. two crossing changes).                              
A link $L$ in $S^3$ is called {\it almost alternating} if
$L$ is not alternating but
admits an almost alternating diagram.
If $D$ is an almost alternating diagram, the
specific crossing to change is called the {\it dealternator} and
we call the other crossings the {\it alternators}.

In this paper, we extend
Theorems \ref{thm:fiber-hopf} and \ref{thm:alt-fiber}
to almost alternating links.
Note that almost alternating diagrams, however, do not always
yield a minimal genus Seifert surface via Seifert's algorithm.
Our first result is as follows:       

\begin{thm} \label{thm:main}
Let $D$ be an almost alternating diagram,
and $R$ a Seifert surface obtained by applying
Seifert's algorithm to $D$. Then, 
$R$ is a fiber surface if and only if
$R$ is connected and desums into a union of Hopf bands.
\end{thm}              

In section 5, we show a stronger version of Theorem \ref{thm:main}
as below,
by using Corollary \ref{cor:algorithm} obtained from the arguments
in the proof of Theorem \ref{thm:main}.

\begin{thm} \label{thm:hopf-plumbing}
Let $R$ be a Seifert surface obtained 
by applying Seifert's algorithm
to an almost alternating diagram.
Then, $R$ is a fiber surface
if and only if $R$ is a Hopf plumbing.
\end{thm}

As a corollary of the proof of Theorem \ref{thm:main},
we obtain
a practical algorithm to determine whether
or not a given almost alternating diagram yields a fiber
surface via Seifert's algorithm.
We use this to prove Theorem 1.5. 
We say that a diagram $D$ is {\it unnested}
if $D$ has no Seifert circle which contains
another circle in both of its complementary region. 
Otherwise we say $D$ is {\it nested}.

\begin{cor}\label{cor:algorithm}
Let $D$ be an
almost alternating diagram and $R$ a Seifert surface
obtained from $D$ by Seifert's algorithm.
Then $R$ is a fiber surface if and only if $R$ 
is connected and desums into a union of
Hopf bands by repeating 
of the following decompositions;\\
(1) a Murasugi decomposition along a nested Seifert circle,\\
(2) a prime decomposition, and\\
(3) Murasugi decompositions of type (A) and (B) in Figure 1.1, 
where each decomposition yields
Seifert surfaces with first Betti numbers smaller than
that of $R$.
\end{cor}

\figbox{Figure 1.1}

\spf
In the proof of Theorems \ref{thm:fiber-hopf} (see \cite[p.533]{G}) 
and \ref{thm:main}, we explicitly show
how we can desum such $R$ into surfaces of smaller first
Betti numbers.
All necessary decompositions are covered in  
the above three.
\fpf

In \cite{Ha}, J. Harer proved that every fiber surface in $S^3$
results from a disk by a sequence of elementary changes
as follows:

(a) plumb on a Hopf band, \par
(b) deplumb a Hopf band, and\par
(c) perform a Dehn twist about a suitable unknotted
curve in the fiber.\par

Then he asked whether changes of either type (b) or (c)
can be omitted, and any fiber surface can be realized only 
using changes of the remaining two types.

So it is worthy presenting the following partial 
affirmative answer as a corollary,
which immediately follows from Theorem 
\ref{thm:hopf-plumbing} and Propositions 5.1 and 5.2.

\begin{cor}\label{cor:only-plumb}
Let $R_1$ and $R_2$ be any fiber surface obtained by applying
Seifert's algorithm to an alternating or almost alternating
diagram.
Then $R_1$ and $R_2$ can be changed into each other by
plumbing and deplumbing Hopf bands.
\end{cor}

We say that a Hopf hand $B$ is {\it positive} (resp. {\it negative})
if the linking number of $\partial B$ is $1$ (resp. $-1$).
By the following fact together with an observation of the
way fiber surfaces deplumb in the proof of Theorem \ref{thm:main}
(see Section 4) and Theorem \ref{thm:hopf-plumbing}
(see Section 5), we have the following Corollary:

\begin{cor}\label{cor:posi-hopf}
Let $D$ be an unnested almost alternating diagram such that the
sign of the dealternator is negative.
Suppose the surface $R$ obtained from $D$ by Seifert's algorithm
is a fiber surface.
Then $R$ is a plumbing of positive Hopf bands.    
\end{cor}

\begin{fact}
Suppose a diagram $D$ is unnested.
Then $D$ is alternating
(resp. almost alternating)
if and only if all the crossings of $D$ have the same sign
(resp. the same sign except exactly one crossing).
\end{fact}

This paper is organized as follows;
Section 2 is for preliminaries.
In Section 3, we give an example for our theorem.
We also show that our theorem can not be extended 
to $2$-almost alternating diagrams, i.e.,
(1) we recall Gabai's example (in \cite{G})
of a $2$-almost alternating diagram for
 a link whose Seifert surface obtained by Seifert's algorithm is
a fiber surface that is not a nontrivial Murasugi sum, and
(2) we give examples of $2$-almost alternating diagrams for knots
whose Seifert surfaces obtained by Seifert's algorithm are fiber
surfaces that are not Hopf plumbing.
In Sections 4 and 5, we prove Theorems~\ref{thm:main} and 
~\ref{thm:hopf-plumbing} respectively.

\section{Preliminaries}

For the definitions of standard terms of sutured manifolds,
see \cite[p.520]{G}.
We say that a sutured manifold $(M,\gamma)$ is a
{\it product sutured manifold} if
$(M,\gamma)$ is homeomorphic to $(R \times I, \partial R \times I)$
with $R_{+}(\gamma)=R \times \{ 1\}, R_{-}(\gamma)=R \times \{ 0\}$,
where $R$ is a compact oriented surface with no closed components
and $I$ is the unit interval $[0,1]$.

The {\it exterior} $E(L)$ of a link $L$ in $S^{3}$
is the closure of $S^{3}-N(L;S^{3})$.
If $R$ is a Seifert surface for $L$, we may assume
$R\cap E(L)$ is homeomorphic to $R$, and often abbreviate 
$R\cap E(L)$ as $R$.

Let $R$ be a Seifert surface for $L$ in $S^3$.
The product sutured manifold
$(M,\gamma)=(R \times I, \partial R \times I)$ 
is called the sutured manifold {\it obtained} from $R$
and the sutured manifold 
$(N,\delta)=(E(L) - {\rm Int}M, \partial E(L) - {\rm Int}\gamma)$
is the {\it complementary} sutured manifold for $R$ 
(or for $(M, \gamma)$).

Note that $R$ is a fiber surface if and only if the 
complementary sutured manifold 
for $R$ is a product sutured manifold.

A {\it product decomposition}
\cite{G}
is a sutured manifold decomposition
$$(M_1,\gamma_1) 
\overset{B}{\longrightarrow} 
(M_2,\gamma_2),$$
where $B$ is a disk properly embedded in $M_1$ such that
$B \cap s(\gamma_1)=$ ($2$ points), 
$M_2 = M_1 - {\rm Int}N(B)$ and that $s(\gamma_2)$ is obtained
by extending $s(\gamma_1) - {\rm Int}N(B)$ in the natural way
(Figure 2.1 (a)).
The disk $B$ is called a {\it product disk}.

Dually, {\it $C$-product decomposition} is the operation
$$(M_1, \gamma_1) 
\overset{E}{\longrightarrow} 
(M_2, \gamma_2),$$
where $E$ is a disk properly embedded in $S^3 - {\rm Int}M_1$ 
such that
$E \cap s(\gamma_1)=$ ($2$ points),
$M_2$ is obtained from $M_1$ by attaching the 2-handle $N(E)$
and that $s(\gamma_2)$ is obtained
by extending $s(\gamma_1) - {\rm Int}N(E)$ in the natural way
(Figure 2.1 (b)).
The disk $E$ is called a {\it $C$-product disk}.

\figbox{Figure 2.1}

\noindent
{\bf Definition.}
Let $R$ be a Seifert surface for a link $L$.
We say that {\it $R$ has a product decomposition} if
there exists a sequence of $C$-product decompositions 
$$(R \times I, \partial R \times I)=(M_0, \gamma_0) 
\overset{E_1}{\longrightarrow}
(M_1, \gamma_1)
\overset{E_2}{\longrightarrow}
\cdots
\overset{E_p}{\longrightarrow}
(M_p, \gamma_p),$$
where the complementary sutured manifold for
$(M_p, \gamma_p)$ is a union of $3$-balls each with a
single suture.

As a criterion to detect a fiber surface, Gabai has shown the
following:

\begin{thm}[\mbox{\cite[Theorem 1.9]{G}}] \label{thm:fiber-pd}
Let $L$ be an oriented link in $S^{3}$, and
$R$ a Seifert surface for $L$. Then,
$L$ is a fibered link with fiber $R$ if and only if
$R$ has a product decomposition.
\end{thm}

We note that in Section 4,
the existence of a $C$-product decomposition 
$(M_0, \gamma_0) 
\overset{E_1}{\longrightarrow}
(M_1, \gamma_1)$ together with the $C$-product disk $E_1$
is important.

\medskip
\noindent
{\bf Definition.}
A surface $R\,(\subset S^{3}$) is a {\it $2n$-Murasugi sum}
of two surfaces $R_{1}$ and $R_{2}$ in $S^{3}$ 
if the following conditions are satisfied;
\begin{enumerate}
\item
$R=R_{1}\underset{\Delta}{\cup}R_{2},$ where $\Delta$ is a $2n$-gon, i.e.,
$\partial \Delta=\mu_{1}\cup\nu_{1}\cup\ldots\cup\mu_{n}\cup\nu_{n}$
(possibly $n=1$), where $\mu_{i} ($resp. $\nu_{i}$) is an arc
properly embedded in $R_{1} ($resp. $R_{2}$).
\item
There exist $3$-balls $B_{1}$ and $B_{2}$ in $S^{3}$ such that:\\
(i) $B_{1}\cup B_{2}=S^{3},\,B_{1}\cap B_{2}=\partial B_{1}=\partial
     B_{2}=S^{2}$ : a $2$-sphere,\\
(ii)$R_{1}\subset B_{1}, R_{2}\subset B_{2}$ and
     $R_{1}\cap S^{2}=R_{2}\cap S^{2}=\Delta.$
\end{enumerate}

The $2$-Murasugi sum is known 
as the connected sum, and
the $4$-Murasugi sum is known 
as the plumbing.

\figbox{Figure 2.2}

Concerning alternating and almost alternating tangles, 
we can confirm the following facts.

\begin{fact}\label{fact:tangle-sum}
Suppose a link diagram $D$ is a tangle sum of 
two tangle diagrams $D_1$ and $D_2$.
If $D$ is alternating, then both $D_1$ and $D_2$ are alternating.
And if $D$ is almost alternating, then one of them, say, $D_1$ is
alternating and $D_2$ is almost alternating.
\end{fact}

\begin{fact}\label{fact:tangle}
By connecting neighboring strands running out of an alternating
(resp. almost alternating ) tangle diagram, we obtain an 
alternating (resp. almost alternating ) link diagram. See Figure 2.3.
\end{fact}

\figbox{Figure 2.3}

Then by these two facts, we can confirm the following propositions.
Let $R$ be a Seifert surface obtained by applying Seifert's algorithm 
to a diagram $D$.

\begin{prop}\label{prop:prime}
If an almost alternating diagram $D$ is a connected sum of
two diagrams, then one of them, say, $D_1$ is alternating and
the other, say, $D_2$ is almost alternating.
The Seifert surface $R$ is a $2$-Murasugi sum of 
$R_1$ and $R_2$, where $R_i$ is obtained from $D_i$.
\end{prop}

\begin{prop}\label{prop:nest}
Suppose that an almost alternating diagram $D$ 
has a nested Seifert circle $C$.
Then, along the disk bounded by $C$, $R$ is 
a Murasugi sum of $R_1$ and $R_2$, where $R_1$ (resp. $R_2$) is
obtained from an alternating (resp. almost alternating) diagram.
\end{prop}

\begin{prop}\label{prop:induction}
Suppose that $R$ desums into two surfaces $R_1$ and $R_2$
as illustrated in Figure 1.1, where the left figures in (A) and (B)
are both almost alternating.
Then $R_i$ $(i=1,2)$ 
is obtained from an alternating or almost alternating diagram. 
\end{prop}

\section{Examples}

In this section, we present some examples.
Example 3.1 is for Theorems \ref{thm:main} and \ref{thm:hopf-plumbing}.
Examples 3.2 and 3.3 show our Theorems \ref{thm:main}
and \ref{thm:hopf-plumbing} can not be extended to 
$2$-almost alternating diagrams.
For the names of knots, refer to Rolfsen's book \cite{R}.

\medskip
\noindent
{\bf Example 3.1.}
Figure 3.1 depicts an almost alternating diagram for the
knot $10_{151}$, 
together with a fiber surface $R$
obtained by Seifert's algorithm.
We can observe that $R$ desums into a union of Hopf bands and 
is a Hopf plumbing.

\figbox{Figure 3.1}

\medskip
\noindent
{\bf Example 3.2.}
Let $R$ be the Seifert surface obtained by applying 
Seifert's algorithm to the oriented pretzel link diagram of 
type $(2,-2,2p)$ as in Figure 3.2, where $p\neq 0$.
$R$ is a fiber surface but does not desum into a union of Hopf 
bands. 

\medskip

We note that this example has been known in \cite{G} 
as a fiber surface for a link 
which does not admit a non-trivial Murasugi sum.

\figbox{Figure 3.2}

\noindent
{\bf Example 3.3.}
Figure 3.3 depicts 2-almost alternating diagrams for the knots
$9_{42}$, $9_{44}$ and $9_{45}$.
By applying Seifert's algorithm to them, we obtain fiber surfaces,
which are not Hopf plumbings.
This can be shown by the following proposition and
direct calculations of genera and the Conway polynomials of
these knots. 

\setcounter{thm}{3}
\begin{prop}[\mbox{\cite[Theorem 3]{MM}}]
If a fibered knot $K$ of genus $2$ can be constructed by plumbing
Hopf bands, then the Conway polynomial $\nabla_K(z)$ of $K$
satisfies the following;
$$
\nabla_K(z) \neq
\begin{cases}
1 + c_{1} z^2 + z^4   & for \ c_1 = 0 \ {\rm mod} \ 4,\cr
1 + c_{1} z^2 - z^4   & for \ c_1 = 2 \ {\rm mod} \ 4. \cr
\end{cases}
$$
\end{prop}

\figbox{Figure 3.3}

\section{Proof of Theorem 1.4}

Since the \lq if\rq \  part is shown by 
Theorem \ref{thm:fiber-sum},
we show the \lq only if\rq \  part.
Let $D$ be an almost alternating diagram for a link $L\,(\subset S^3)$
on the {\it level $2$-sphere} $S^2$ and
let $R$  be a Seifert surface obtained
by applying Seifert's algorithm to $D$.

Note that if a diagram $D$ is unnested, then Seifert's algorithm uniquely
yields a Seifert surface. We say that a  Seifert surface $R$
is {\it flat} if $R$ is obtained from an unnested diagram and thus lies
in $S^2$ except in the neighborhood of each crossing. 

Suppose that $R$ is a fiber surface. 
Since any fiber surface is connected,
we can assume $D$ 
is connected.

Suppose $D$ is nested.
Then, by Theorem  \ref{thm:fiber-sum},
$R$ desums into fiber surfaces $R_1$ and $R_2$.
Moreover, by Proposition \ref{prop:nest}, 
one of them, say, $R_1$ is obtained from an alternating 
diagram and $R_2$ from an almost alternating diagram.
By Theorem \ref{thm:fiber-hopf}, 
$R_1$ desums into a union of Hopf bands.
Therefore, we may assume that $D$ is unnested.

Similarly, by Proposition \ref{prop:prime}, 
we may assume that $D$ is prime, and in particular, reduced.

Now we prove the theorem by induction
on the first Betti number $\beta_{1}$ of $R$, 
where $R$ is a fiber surface obtained by applying 
Seifert's algorithm to a connected unnested prime almost 
alternating diagram $D$. 
If $\beta_{1}=1$, then $R$ is an unknotted annulus
and $D$ has $n$ crossings which are of the same sign
except exactly one crossing.
Note that $R$ is a fiber surface if and only if $n=4$,
in which case $R$ is a Hopf band.
Hence we have the conclusion.

Then we assume that the theorem holds when $\beta_1(R) < k$
and prove the theorem for $R$ with $1\le\beta_1(R) = k$.

The main method of 
the proof is to examine the $C$-product disk for the
sutured manifold obtained from $R$ and grasp a local picture
where we can desum $R$ into
surfaces $R_1$ and $R_2$ obtained by the algorithm 
with smaller first Betti numbers.
In each case, 
it is easy to confirm that $D_{i}\,(i=1,2)$ is an alternating 
or almost alternating diagram, that is, they satisfy
the assumption of the induction 
(see Corollary \ref{cor:algorithm} and Proposition \ref{prop:induction}).

Let $(M, \gamma)$ be the sutured manifold obtained from $R$.
We identify $s(\gamma)$ as $L$.
Let $E$ be a $C$-product disk
for $(M, \gamma)$, i.e., 
$E$ is properly embedded in 
$S^3 - {\rm Int}M$ so that $E \cap L =$ ($2$ points).
We may suppose that $E$ is non-boundary-parallel, and 
assume that $|E\cap S^2|$ is minimal
among all such disks.
Further, we may assume by isotopy that $\partial E \cap L$ 
occurs only in small neighborhoods of the crossings of $D$.
Similarly, we can assume that $\partial E\cap S^{2}$ occurs only
in small neighborhoods of the crossings.
For convenience, we say that
{\it $\partial E\cap L$ and $\partial E\cap S^{2}$ occur at the crossings.} 

{\bf Case A.} 
$E\cap S^{2}=\emptyset.$

If $\partial E\cap L$ occurs at one crossing, then
$E$ is boundary parallel, a contradiction. 
Thus, we suppose that $\partial E\cap L$ occurs at two crossings
(see Figure 4.1).
If both crossings are alternators,
we see that $R$ is a plumbing of flat surfaces, 
one of which is obtained from an unnested almost 
alternating diagram and has first Betti number 
smaller than $k$.
If one crossing is the dealternator,
we also see that $R$ is a plumbing of surfaces,
one of which is compressible and hence not a fiber surface,
a contradiction to Theorem \ref{thm:fiber-sum}.

\figbox{Figure 4.1}

{\bf Case B.} 
$E\cap S^{2}\neq\emptyset.$

Label the crossings with
$\fr0, \fr1, \ldots , \fr{w-1}$
so that the dealternator has $\fr0$.
By standard innermost circle argument, we may assume,
by the minimality of $|E \cap S^2|$, that 
$E\cap S^{2}$ consists of arcs.
Let $\alpha$ be an arc of $E\cap S^{2}$.
By assumption, each endpoint of $\alpha$ lies in a neighborhood of a 
crossing and hence is accordingly labeled.
Then the {\it label} of $\alpha$ is a pair $(\fr{i}, \fr{j})$
of the labels of $\partial \alpha$.
The two points of $\partial E \cap L$
are also labeled according to the crossings
at which $\partial E \cap L$ occurs.

\begin{lem}
For any arc $\alpha$ of $E \cap S^2$ with label
$(\fr{i}, \fr{j})$, we have $i\neq j$.
\end{lem}

\spf
If both of the endpoints of $\alpha$ occur at 
the same crossing $\fr{i}$,
we can observe that one of the two cases in Figure 4.2 occurs.
In Figure 4.2 (a), $D$ is non-prime.
In Figure 4.2 (b), there exists an arc $\alpha '$ of $E\cap S^2$ in
 $S^2 - {\rm Int}M$ such that
the endpoints of $\alpha '$ occur at the same crossing $\fr{i}$, and that
$\alpha '$ cuts off a disk $H$ from $S^2-\text{Int}M$
with $\text{Int}H\cap(E\cap S^2)=\emptyset$.
We can surgery $E$ along $H$
so that we obtain two disks $E_1,\,E_2$
properly embedded in $S^3 - {\rm Int}M$.
Since both endpoints of $\alpha '$ are in $R_{+}(\gamma)$
(or $R_{-}(\gamma)$), one of them, say $E_1$, intersects $L$ twice.
Since $E$ is non-boundary-parallel, so is $E_1$ or $E_2$.
If $E_2$ is, then it yields 
a compressing disk for $R$, 
a contradiction.
Hence $E_1$ is a non-boundary-parallel $C$-product disk
with $|E_1\cap S^2|<|E\cap S^2|$, a contradiction.
\fpf

\figbox{Figure 4.2}

We look at an outermost disk $F \subset E$
(i.e., $F$ is the closure of a component of $E - S^2$
  such that $F \cap S^2$ is connected).

\begin{lem}
Let $\alpha$ be an outermost arc of $E\cap S^2$ 
with label $(\fr{i}, \fr{j})$, 
cutting an outermost disk $F$ off $E$. 
Then we may assume that $i\neq j$ and that 
$i$ or $j=0$ if 
$F\cap L= \emptyset$ or (a point).
\end{lem}

\spf By Lemma 4.1, we have $i \neq j$.
Suppose $i \neq 0$ and $j \neq 0$.
If $|F\cap L|=0$, $R$ is non-prime, a contradiction
(Figure 4.3 (a)).
If $|F\cap L|=1$, then either $|E \cap S^2|$ is not minimal,
or $R$ is a plumbing (Figures 4.3 (b) and (c)).
\fpf

\figbox{Figure 4.3}

Concerning outermost disks,
we have two cases.

{\bf Case B-1.}
{\it There exists an outermost disk $F$ with $F \cap L =$ (a point).}\\
Let $\alpha$ be the arc $F \cap S^2\,(\subset E)$.
By Lemma 4.2, we assume the label of $\alpha$ is
$(\fr{0}, \fr{j})$, where $j \neq 0$.
Let $\fr{k}$ be the label of the point of $\partial E \cap L$ on $F$.
Then we have three cases;
{\bf Subcase 1}: $k = 0$,
{\bf Subcase 2}: $k = j$, and
{\bf Subcase 3}: $k \neq 0$ and $k \neq j$.
In Subcases 1 and 2, $D$ is non-prime (Figure 4.4 (a)).
In Subcase 3, $R$ is a plumbing or we can isotope $E$
so that the outermost disk of Case B-1 is replaced
by an outermost disk of Case B-2 (Figure 4.4 (b)).

\figbox{Figure 4.4}

\begin{lem}
We may assume there exists no outermost disk of Case B-1.
\end{lem}

\spf
If the latter situation of Subcase 3 
above occurs,
we can view the above isotopy of $E$ as sliding a point of
$\partial E \cap L$ out of $F$. Hence by repeating the above 
isotopies at most twice, 
we may eliminate outermost disks of Case B-1.
\fpf

{\bf Case B-2.}
{\it There exists an outermost disk $F$ with
  $F \cap L = \emptyset$.}\\
By Lemma 4.2,
we may assume $\alpha = F \cap S^2$
appears as in Figure 4.5.
We note that outermost disks of this kind are typically found in
the complementary sutured manifold for the fiber surface in Figure 3.2,
which is obtained from a $2$-almost alternating diagram. 
The rest of the proof really depends on the
almost-alternatingness of $D$.

\figbox{Figure 4.5}

\begin{lem}
For any arc $\beta$ of $\partial E - (E \cap S^2)$,
if $\beta \cap L = \emptyset$, then
the endpoints of $\beta$ have different labels.
\end{lem}

\spf
Suppose the 
two endpoints of $\beta$ have the same 
label. Then $\beta$ appears as in Figure 4.6 and we can 
isotope $E$ to a $C$-product disk $E'$ such that 
$|E'\cap S^2|=|E\cap S^2|-1$, a contradiction.
\fpf

\figbox{Figure 4.6}

\begin{lem}
Suppose $E$ locally appears as in Figure 4.7 (a), 
i.e., $\fr{i},\ \fr{j}$ and $\fr{k}$ are the labels of
points of
$(\partial E\cap S^2) \cup (\partial E \cap L)$
sequential in $\partial E$
such that the former two 
points are connected by an outermost arc of $E\cap S^2$
and the last is a point of $\partial E \cap S^2$.
Then $i, j, k$ are mutually different.
\end{lem}

\spf
By Lemma 4.2, we have $i\neq j$ and $i$ or $j=0$. 
Suppose $i=k$. Then we can find a compressing 
disk for $R$ in Figure 4.7 (b), a contradiction. 
By Lemma 4.4, we have $j\neq k$.
\fpf

\figbox{Figure 4.7}

\begin{lem}
We may assume that the following situation never occurs;
The disk $E$ locally appears as in Figure 4.8 (a), i.e.,
$\fr{i},\,\fr{j},\,\fr{k}$ and $\fr{l}$ are the
labels of points of $(\partial E\cap S^2) \cup (\partial E \cap L)$
sequential in $\partial E$
such that the former two points and the latter two
are respectively connected by outermost arcs 
$\alpha_1$ and $\alpha_2$ of $E\cap S^2$.
\end{lem}

\spf
By Lemmas 4.2, 4.4 and 4.5, we may assume
that $i=l=0$ and $j\neq k$.
Then we obtain the conclusion, since in Figure 4.8 (b)
$\alpha_2$ can not coexist with the arc of $\partial E-S^2$ 
connecting $\fr{k}$ and $\fr{l}$.
\fpf

\figbox{Figure 4.8}

\begin{lem}
Suppose $E$ locally appears as in Figure 4.9 (a), i.e.,
$\fr{i}, \fr{j}$ and $\fr{k}$ are the labels of points of
$(\partial E \cap S^2) \cup (\partial E \cap L)$
sequential in $\partial E$
such that the former two points are connected
by an outermost arc of $E \cap S^2$ and
the third point is of $\partial E \cap L$.
Then we may assume $i, j, k$ are mutually different.
\end{lem}

\spf
By Lemma 4.2, we have $i \neq j$, and $i$ or $j =0$.
If $k=i$, then $R$ is compressible, a contradiction
(see Figure 4.9 (b)).
If $k=j$, we can reduce $|E\cap S^2|$ by isotopy, a contradiction 
(see Figure 4.9 (c)).
\fpf

\figbox{Figure 4.9}

\begin{lem}
Let $\fr{l}$ be the label of the point $x$ of
$\partial E\cap L$.
Suppose that the two points adjacent to $x$ in $\partial E$
are points of $\partial E\cap S^2$.
Then the two adjacent points
do not have the same label except for
the case where they are both $\fr{l}$.
\end{lem}

\spf
Suppose the two points have the same label $\fr{i} (\neq \fr{l})$.
By Lemma 4.1, we may assume that they are not connected
by an arc of $E\cap S^2$.
Then we can find a $C$-product disk $E'$ in Figure 4.10
such that  $|E'\cap S^2|=0$, a contradiction.
\fpf

\figbox{Figure 4.10}

Similarly we have the following lemma.

\begin{lem}
Suppose $E$ locally appears as in Figure 4.11(a),
i.e., $\fr{i}, \fr{j}, \fr{k}$ and $\fr{l}$ are the 
labels of points of
$(\partial E \cap S^2) \cup (\partial E \cap L)$
sequential in $\partial E$ such that 
the former two points are connected
by an outermost arc of $E \cap S^2$,
the third point is of $\partial E \cap L$ and that 
the fourth is of $\partial E \cap S^2$.
If $k \neq l$, then $i, j, k$ and $l$ are mutually different.
\end{lem}

\spf
By Lemma 4.7,
we may assume $i, j, k$ are mutually different.
Then by Lemma 4.8, we have $l \neq j$.
Suppose $l \neq k$ and $l = i$.
Then by Lemma 4.2, 
$\fr{i} = \fr{l} = \fr{0}$ (Figure 4.11 (b)) or $\fr{j}=\fr{0}$ (c).
In either case,  we can find a $C$-product disk $E'$
such that $|E'\cap S^2|=0$, a contradiction.
\fpf

\figbox{Figure 4.11}

\begin{lem}
We may assume the following situation never occurs;
The disk $E$ locally appears as in Figure 4.12 (a),
i.e., $\fr{i}, \fr{j}, \fr{k}, \fr{l}$ and $\fr{m}$
are the labels of points of
$(\partial E \cap S^2) \cup (\partial E \cap L)$
sequential in $\partial E$
such that 
the first two points 
and the last two points
are respectively connected
by an outermost arc of $E \cap S^2$,
and that the third point is of $\partial E \cap L$.
\end{lem}

\spf
By Lemma 4.7, $k, l$ and $m$ are mutually different
and hence by Lemma 4.9, $i, j, k$ and $l$ are mutually different.
By Lemma 4.2, $i$ or $j=0$ and $l$ or $m=0$, and hence $m=0$, 
and by symmetry, we have $i=0$.
Then $\fr{i}, \fr{j}, \fr{k}, \fr{l}$ and $\fr{m}$ appear as in
Figure 4.12 (b), where
we can find a $C$-product disk $E'$ 
such that $|E \cap S^2| > |E'\cap S^2|=1$, a contradiction.
We note that $E'\cap L$ occurs at 
$\fr{k}$ and $\fr{l}$.
\fpf

\figbox{Figure 4.12}

An arc $\varepsilon$ of $E \cap S^2$ is 
said to be {\it of level $2$} if 
it is not outermost and, 
for one component $E_1$ of $E - \varepsilon$,
$E_1 \cap S^2$ is a union of outermost arcs in $E \cap S^2$.
Suppose there is no arc of level $2$.
Then by Lemmas 4.3, 4.6 and 4.10,
we see that $E \cap S^2$ consists of only one arc $\alpha$
such that one component of $E - \alpha$
contains the two points of $\partial E \cap L$.
Let $(\fr{0},\fr{j})$ be the label of $\alpha$,
and let $\fr{k}$ and $\fr{l}$ be the labels of the two points of
$\partial E\cap L$, where $\fr{0}, \fr{j}, \fr{k}$ and $\fr{l}$
appear in this order in $\partial E$.
If $l=k$, then we can isotope $E$ so that $E \cap L = \emptyset$
and we have a compressing disk for $R$, for
$E$ is not boundary parallel, a contradiction. 
Hence by Lemma 4.7, we can assume
$j, k, l, 0$ are mutually different.
In this case, $R$ desums into three surfaces
$R_1, R_2$ and $R_3$
obtained by applying Seifert's algorithm to
the almost alternating diagrams $D_1, D_2$
and $D_3$ respectively (Figure 4.13).

\figbox{Figure 4.13}

Hence we assume there is an arc of level $2$. 
Then by Lemmas 4.3, 4.6 and 4.10, we see that 
there exists an arc $\varepsilon$ of level $2$ such that
one disk $E_{1}$ cut by $\varepsilon$ off $E$
contains one outermost arc of $E\cap S^2$ and 
satisfies one of the following conditions;\\
(*) $E_1\cap L=\emptyset$, \\
(**) $E_1\cap L=$ a point. 

If $E_{1}$ satisfies (*), by Lemmas 4.1 and 4.5, 
all four labels of
points of $E_1\cap S^2$ are mutually different.
Then, we can see that $D$ is non-prime or $R$ is 
a plumbing (Figure 4.14).

\figbox{Figure 4.14}

Thus we have:

\begin{lem}
We may assume that there is 
no arc of level $2$ which cuts a disk $E_{1}$ off E
such that $E_{1}$ contains only one (outermost) 
arc of $E\cap S^2$ and that
$E_{1}\cap L=\emptyset$.
\end{lem}

In what follows, we assume that there exists 
an arc $\varepsilon$ of level $2$
which cuts off $E$ a disk $E_{1}$ containing
one outermost arc of $E\cap S^2$ and 
satisfying (**).

By Lemma 4.3, we may suppose that 
$E_1$ appears as in Figure 4.15 (a) with labels
$\fr{i},\,\fr{j},\,\fr{k},\,\fr{l}$ and $\fr{m}$.

\begin{lem}
All five labels in $E_{1}$ are mutually different.
\end{lem}

\spf
By Lemma 4.5, $i, j, k$ are mutually different.
By Lemma 4.7, $l \neq k$ and $l \neq j$.
We see $l \neq i$, for if not, $R$ appears as in
Figure 4.15 (b) or (c), and in either case, $R$ is compressible, 
a contradiction. 
Now we have seen that $i, j, k, l$ are mutually different.
Next suppose $m = l$. Then $R$ appears as in Figure 4.15 (d) or (e).
In Figure 4.15 (d), $R$ is a plumbing 
or we can isotope $E$ to 
reduce $|E \cap S^2|$. In Figure 4.15 (e), 
$R$ is a Murasugi sum
or we can isotope $R$ so that $D$ becomes an alternating diagram
and the result follows from Theorem \ref{thm:fiber-hopf}.
Hence we can assume $m \neq l$ and by Lemma 4.9, we see that
$j, k, l, m$ are mutually different and by Lemma 4.1, $m \neq i$.
\fpf

\figbox{Figure 4.15}

\begin{lem} We may assume
$\fr{j}=\fr{0}$.
\end{lem}

\spf
If not, $\fr{k}=\fr{0}$ by Lemma 4.2.
Then $R$ is a 6-Murasugi sum as in Figure 4.16.
\fpf

\figbox{Figure 4.16}

\begin{lem}
Let $\varepsilon$ and $E_1$ be as above. 
Then there is no arc $\varepsilon '$ of $E \cap S^2$ as in Figure 4.17
which cuts a disk $E_2$ off $E$ with the following conditions:
\begin{enumerate}
\item
$E_{1}\subset E_{2}$,
\item
$({\rm Int }E_{2}-E_{1})\cap(E\cap S^2)=\emptyset$,
\item
$E_{2}\cap L=E_{1}\cap L= (1\,{\rm point})$.
\end{enumerate}
\end{lem}

\figbox{Figure 4.17}

\spf
By Lemmas 4.12 and 4.13, 
we may assume that $E_1$ appears as in Figure 4.18. 
Recall that $R$ is flat. 
Suppose that we have a disk $E_2$ as in Figure 4.17.
Then the arc $\varepsilon'$ lies in some
region of $S^2 - N(R)$. Hence, considering the orientation
of $R$, we see that one of the following occurs;\\
(1) The point \textcircled{\small 1} is bounded by the same
Seifert circle as one of the points
\textcircled{\small 3} and \textcircled{\small 5}, \\
(2) The point \textcircled{\small 6} is bounded by the same
Seifert circle as one of the points
\textcircled{\small 2} and \textcircled{\small 4}.\\
In each case, we can find a $C$-product disk $E'$ 
such that $|E\cap S^2| > |E'\cap S^2| = 0 \ {\rm or\ } 1$, a contradiction.
\fpf

\figbox{Figure 4.18}

\begin{lem}\label{lem:last}
Let $E_1$ be as above.
Then the following situation never occurs;
The disk $E$ locally appears as in Figure 4.19,
i.e., there is an outermost disk $F$ such that 
$\partial E - (E_1 \cup F)$
has a component $\beta$ which contains no point of 
$(\partial E\cap S^2) \cup (\partial E\cap L)$.
\end{lem}

\figbox{Figure 4.19}

\spf
Suppose there exists such a disk $F$.
Let $\alpha$ be an arc in $E\cap S^2$ which cuts $F$ off $E$,
and $(\fr{s},\fr{t})$ the label of $\alpha$
where $\fr{s}$ is the label of an endpoint of $\beta$.
First we examine 
the case where $E$ appears as in Figure 4.19 (a).
If $s=i$ or $0$, we
can find a $C$-product disk $E'$ such that $|E'\cap S^2|=1$, 
a contradiction (Figure 4.20 (a)).
By Lemma 4.2, we have $s=0$ or $t=0$, and hence $t=0$.
If $s = k$ or $l$, we can find a 
$C$-product disk $E'$ such that $|E' \cap S^2| = 0$, a contradiction
(see Figure 4.20 (b)).
By Lemma 4.4, we have $s \neq m$. 
Then we see that $R$ locally appears 
as in Figure 4.20 (c).
It is impossible that $\partial E$ runs 
toward the dealternator 
$\fr{0}(=\fr{t})$ after passing through $\fr{s}$
because of the orientation of $R$. 

Second, we examine the case where
$E$ locally appears as in Figure 4.19 (b).
We can do this by the similar way to 
in the previous case.
By Lemma 4.4, $s \neq i$.
If $s = 0$ or $k$, we can find a $C$-product disk $E'$ 
such that $|E'\cap S^2|=0$, a contradiction.
By Lemma 4.2, we have $s=0$ or $t=0$, and hence $t=0$.
If $s = l$ or $m$, we can find a $C$-product disk $E'$ 
such that $|E'\cap S^2|=1$, a contradiction. 
Then we see that it is impossible that $\partial E$ 
runs toward the dealternator $\fr{0}(=\fr{t})$ before 
passing through $\fr{s}$. 
See Figure 4.21.
\fpf

\figbox{Figure 4.20}

\figbox{Figure 4.21}

Let $E_1' = E - E_1$. 
Then $E_1' \cap L$ is exactly one point, say, $x$.
By Lemma 4.3, 
$E_1' \cap (E \cap S^2) \neq \emptyset$.
By Lemmas 4.6 and 4.11,
any arc of $E_1' \cap (E \cap S^2)$ which does not separate
$\varepsilon$ and $x$ is outermost in $E_1'$.
By Lemma 4.15, 
at least one of $E_1' \cap (E \cap S^2)$
separates $\varepsilon$ and $x$.
Among such separating arcs, let $\alpha$ be the one
closest to $\varepsilon$.
Then by Lemma 4.15 again,
the subdisk of $E$ between $\varepsilon$ and $\alpha$
contains no arc of $E_1' \cap (E \cap S^2)$.
However, this contradicts Lemma 4.14.
This completes the proof.
\epf{\ref{thm:main}}

\section{Proof of Theorem 1.5}

In this section, we prove Theorem \ref{thm:hopf-plumbing}.
Recall that a Seifert surface $R$ obtained by Seifert's algorithm
is a union of {\it Seifert disks} and {\it Seifert bands}.

\medskip
\noindent
{\bf Definition.} 
Let $R$ be a Seifert surface obtained by Seifert's algorithm.
We say that two Seifert bands $B_1$ and $B_2$ of $R$ are {\it parallel}
if they connect the same two Seifert disks.

\medskip
The following is a case where we can deplumb a Hopf
band from a fiber surface:

\begin{prop} \label{prop:para-deplumb}
Let $R$ be a fiber surface obtained by Seifert's algorithm.
Suppose $R$ has a pair of parallel bands $B_1$ and $B_2$.
Then, we can deplumb a Hopf band from $R$.
Moreover, we have the following;\\
(1) the parallel bands are of the same sign, and\\
(2) for each $i = 1, 2$, 
we can cut the band $B_i$ by deplumbing a Hopf band from $R$,
i.e., $R$ is a plumbing of $R-B_{i}$  and a Hopf band.
\end{prop}

\noindent
\spf
We denote by $L$ the link $\partial R$.
We may assume that
the Seifert circles, say, $C_1$ and $C_2$
connected by $B_1$ and $B_2$ bound 
mutually disjoint Seifert disks 
on the level 2-sphere $S^2$.

First, suppose the pair of parallel bands 
are of the same sign. 
We may assume they appear as in Figure 5.1 (a). 
We explicitly show that $R$ is a plumbing of a Hopf band 
and the surface $R-B_{i}$.
Move $L$ by isotopy as in Figure 5.1 (a) 
and let $R'$ be the surface as depicted.
Apparently the Euler characteristic $\chi(R)$ is equal to
$\chi(R')$.
Hence by the uniqueness of fiber surfaces,
we see that $R$ is isotopic to $R'$.
Now we can deplumb a Hopf band from $R'$ 
as in Figure 5.1 (b).
Then by retracing the above isotopy, we obtain the conclusion.

Next suppose that the pair of parallel bands are of the 
opposite signs, i.e., that the twisting of $B_1$ is opposite.
Then by the isotopy as implied by Figure 5.1 (a), 
we can find a compressing disk for $R'$, which contradicts
the fact that fiber surfaces are of minimal genus and hence
incompressible.
\fpf

\figbox{Figure 5.1}

The following proposition assures that
if a diagram $D$ has a Seifert circle $C$ which contains an alternating
tangle diagram, then any Seifert surface obtained by applying Seifert's
algorithm to $D$ has parallel bands.

\begin{prop} \label{prop:alt-para}
Suppose a Seifert surface $R$ obtained from an alternating 
diagram $D$ is a fiber surface.
Then $R$ has parallel bands.
Moreover, if $D$ is reduced, then for any band $B$ of $R$,
there is a band $B'$ of $R$ which is parallel to $B$.
\end{prop}

\spf
By untwisting $R$ by isotopy if necessary, 
we may assume that
$D$ is reduced.
Moreover, we may assume that $D$ is unnested, because
    (1) by desumming along nested Seifert circles, we can 
decompose $R$ into fiber surfaces obtained from 
unnested alternating diagrams, and
    (2) if one of the decomposed surfaces has parallel bands,
then so does $R$.
Suppose a fiber surface $R$ for a link $L$
is obtained from a reduced unnested alternating diagram $D$.
Then by \cite{Mu1960}
(or \cite[Proposition 13.25]{BZ}), $L$ is a connected sum of
$(2, n)$-torus knots or links. 
Moreover the arguments in \cite{Mu1960} shows that $D$ is 
the \lq standard' alternating diagram of a connected sum of
$(2, n)$-torus knots or links. 
Hence we obtain the conclusion.
\fpf

\medskip 
\noindent
{\it Proof of Theorem \ref{thm:hopf-plumbing}.}
The \lq if' part follows from Theorem \ref{thm:fiber-sum}.
We show the \lq only if' part, using 
Corollary \ref{cor:algorithm},
by induction on the first Betti
number $\beta_1$ of $R$. 
If $\beta_1(R) = 1$, $R$ is a Hopf band, and hence 
the theorem holds. 
Assume the theorem holds for such surfaces
with $\beta_1 < k$, and let
$R$ be a Seifert surface with $\beta_1(R) = k$
obtained from an almost alternating diagram $D$.
By untwisting $R$ if necessary, we may assume that
$D$ is reduced.
By Corollary \ref{cor:algorithm}, we know how $R$
decomposes into Hopf bands. 
Hence by the following four lemmas,
we will see that we can deplumb a Hopf band from $R$,
in such a way that by deplumbing a Hopf band, 
we cut a band of $R$ corresponding to an alternator.
Therefore the deplumbed surface satisfies the
assumption of induction so that 
we see that
$R$ is a Hopf plumbing. 
\epf{\ref{thm:hopf-plumbing}}

\begin{lem} \label{lem:claim1}
If $R$ desums along a nested Seifert circle, then
we can cut a band of $R$ by deplumbing a Hopf band from $R$.
\end{lem}

\spf
Suppose $D$ is nested, i.e., there exists a Seifert circle $C$
which contains another Seifert circle in both of its
complementary regions in $S^2$.
Then $R$ desums along $C$ into two surfaces, 
say, $R_1$ and $R_2$ such that 
$R_1$ is obtained from an alternating diagram
and $R_2$ from an almost alternating diagram
(cf. Proposition \ref{prop:nest}).
Note that by Theorem \ref{thm:fiber-sum}, 
both $R_1$ and $R_2$ are fibers.
By Proposition \ref{prop:alt-para}, we see that $R_1$ has
parallel bands and hence so does $R$.
Then by Proposition \ref{prop:para-deplumb}, we can cut
a band of $R$ by deplumbing a Hopf band from $R$.
\fpf

\begin{lem}\label{lem:claim2}
If $R$ is a connected sum, then 
we can cut a band of $R$ by deplumbing a Hopf band from $R$.
\end{lem}

\spf
Let $R$ be a connected sum of $R_1$ and $R_2$, where
$R_1$ is obtained from an alternating diagram and $R_2$
from an almost alternating diagram by Proposition \ref{prop:prime}.
Then by Theorem \ref{thm:fiber-sum} and 
Proposition \ref{prop:alt-para}, 
$R_1$ has parallel bands, 
which are also parallel in $R$, and hence,
by Proposition \ref{prop:para-deplumb},
we can cut a band of $R$ by deplumbing a Hopf band from $R$.
\fpf

\begin{lem}\label{lem:claim3}
If $R$ admits a decomposition of type (A), then
we can cut a band of $R$ by deplumbing a Hopf band from $R$.
\end{lem}

\spf
Suppose $R$ admits a decomposition of type (A).
Then we can deform $R$ to $R'$ by isotopy 
as depicted in Figure 5.2 (a),
from which we can desum a fiber surface $R_1$ in Figure 5.2 (b).
We can confirm that $R_1$ 
is obtained from an alternating diagram using Fact \ref{fact:tangle}.
By Proposition \ref{prop:alt-para}, $R_1$ has parallel bands.
Though $R'$ itself is not a surface obtained by Seifert's algorithm,
we can apply the argument in the proof of 
Proposition \ref{prop:para-deplumb}, by regarding the inside of
the dotted circle in Figure 5.2 (a) as a black box.
Hence we can cut a band of $R'$ 
(which is a band in the image of $R_1$ in $R'$)
by deplumbing a Hopf band from $R'$. 
This corresponds to cutting a band of $R$ 
by deplumbing a Hopf band from $R$.
Note that we can confirm that the surface obtained 
from $R$ by this cutting the band satisfies the
assumption of induction.
\fpf

\figbox{Figure 5.2}

\begin{lem}\label{lem:claim4}
If $R$ admits a decomposition of type (B), then
we can cut a band of $R$ by deplumbing a Hopf band from $R$.
\end{lem}

\spf
According to whether the crossing visible in Figure 1.1 is
an alternator or the dealternator, we have two cases.
Let us call the former a decomposition of type (B1)
and the latter of type (B2).
Suppose that $R$ admits a decomposition of type (B1).
Then by the same way as in the proof of
Lemma \ref{lem:claim2}, 
we can cut a band of $R$ by deplumbing a Hopf band from $R$.

Now assume $R$ does not admit a decomposition of type (B1).
Then $R$ deplumbs into $R_1$ and $R_2$, which 
are both obtained from almost alternating diagrams 
(see Proposition \ref{prop:induction}).
If $R_1$ or $R_2$ admits a decomposition
of type (A), then we see, 
by the uniqueness of fiber surfaces,
that $R$ also admits a decomposition of type (A), and
the claim follows from Lemma \ref{lem:claim3}.
Hence we assume that neither $R_1$ nor $R_2$ admits a 
decomposition of type (A). 
Inductively, if we can do a decomposition of type (A) or (B1)
in the process of desumming $R$ into a union of Hopf bands, then
we see that $R$ also admits a decomposition of type (A) or (B1).
So we assume that
$R$ desums into a union of Hopf bands using decompositions
of type (B2) alone. Then by another inductive argument, we see that
$R$ is a pretzel surface of type
$(1, -3, \ldots, -3)$ or $(-1, 3, \ldots, 3)$.
In this case, obviously
we can cut a band of $R$ by deplumbing a Hopf band from $R$.
\fpf

{\bf Acknowledgment.}
The authors would like to thank Professor Taizo Kanenobu, 
Professor Tsuyoshi Kobayashi, Professor Yasutaka Nakanishi 
and Professor Makoto Sakuma for their comments.
Part of this work was carried out while the first 
author was visiting at University of California, Davis. 
He would like to express thanks to Professor 
Abigail Thompson and the department for their hospitality.

\footnotesize{

}

Hiroshi Goda \hspace{3mm}
goda@math.kobe-u.ac.jp\\
Graduate School of Science and Technology, 
Kobe University, \\
Rokkodai 1-1, Nada, Kobe, 
657-8501, Japan \\

Mikami Hirasawa \hspace{3mm}
hirasawa@math.sci.osaka-u.ac.jp\\
Department of Mathematics,
Graduate School of Science,\\
Osaka University,
Machikaneyama 1-1, Toyonaka, 560-0043, Japan\\

Ryosuke Yamamoto \hspace{3mm}
ryosuke@math.kobe-u.ac.jp\\
Graduate School of Science and Technology, 
Kobe University, \\
Rokkodai 1-1, Nada, Kobe, 
657-8501, Japan


\begin{thebibliography}{99}

\bibitem{A} C. Adams and et al:
{\it Almost alternating links},
Topology and its Appl., 46 (1992), 151-165.

\bibitem{BZ} G. Burde and H. Zieschang:
{\rm Knots}, Studies in Math. 5, Walter de Gruyter, 1985.

\bibitem{C} R. Crowell, {\it Genus of alternating link types},
 Ann. of Math., 69 (1959), 258-275.

\bibitem{G} D. Gabai:
{\it Detecting fibred links in $S^{3}$},
Comment. Math. Helv., 61 (1986), 519-555.

\bibitem{Ha} J. Harer:
{\it How to construct all fibered knots and links},
Topology, 21 (1982), 263-280.

\bibitem{MM} P. Melvin and H. Morton:
{\it Fibred knots of genus 2 formed by
plumbing Hopf bands},
J. London Math. Soc., (2) 34 (1986), 159-168.

\bibitem{Mu1958} K.\ Murasugi, 
{\it On the genus of the alternating knot, I, II},
 J. Math. Soc. of Japan, 10 (1958), 94-105, 235-248.

\bibitem{Mu1960} K. Murasugi:
{\it On alternating knots},
Osaka Math. J., 12 (1960), 277-303.

\bibitem{Mu} K. Murasugi:
{\it On a certain subgroup of an alternating link},
Amer. J. Math., 85 (1963), 544-550.


\bibitem{R} D. Rolfsen:
{\rm Knots and Links}, Math. Lecture Series 7,
Publish or Perish Inc., Berkeley, 1976.

\bibitem{St} J. Stallings:
{\it On fibering certain $3$-manifolds},
In: Topology of $3$-manifolds and related topics (Georgia, 1961)
Prentice Hall (1962), 95-100.
\bigskip

\end{thebibliography}
\end{document}